\documentclass[1p]{elsarticle}

\journal{Number Theory}
\usepackage{amsmath,amsthm,amssymb}
\usepackage{amsmath,amsfonts,amssymb,amstext,amsthm,amscd,mathrsfs,amsbsy}

\def\noi{\noindent}
\newtheorem{thm}{Theorem}
\newtheorem{lemma}[thm]{Lemma}
\newtheorem{prop}[thm]{Proposition}

\newtheorem{definition}{Definition}

\newdefinition{remark}{Remark}
\newproof{pf}{Proof}
\newproof{p1}{Proof of theorem \ref{thm1}}
\newproof{p2}{Proof of theorem \ref{thm2}}

\theoremstyle{definition}

\theoremstyle{plain}

\newtheorem{coro}{Corollary}[section]

\newcommand{\la}{\langle}
\newcommand{\ra}{\rangle}

\newcommand{\ZZ}{\mathbb{Z}}
\newcommand{\QQ}{\mathbb{Q}}
\newcommand{\RR}{\mathbb{R}}
\newcommand{\OO}{\mathfrak{o}}
\newcommand{\HH}{\mathscr{H}}
\newcommand{\CC}{\mathbb{C}}
\newcommand{\NN}{\mathbb{N}}
\newcommand{\KK}{\mathrm{K}}

\newcommand{\MM}{\mathscr{M}}
\newcommand{\M}{\mathcal{M}}

\newcommand{\R}{\mathrm{R}}
\newcommand{\D}{\mathrm{D}}

\newcommand{\E}{\mathrm{E}}

\newcommand{\zz}{\mathbf{z}}
\newcommand{\xx}{\mathbf{x}}
\newcommand{\yy}{\mathbf{y}}
\newcommand{\X}{\mathrm{X}}
\newcommand{\Y}{\mathrm{Y}}

\newcommand{\T}{\mathrm{T}}
\newcommand{\XX}{\mathcal{H}}

\newcommand{\fa}{\mathfrak{a}}
\newcommand{\fb}{\mathfrak{b}}

\newcommand{\q}{q}

\newcommand{\cu}{\lambda}

	\newcommand{\PSL}{PSL(2,\mathbb{R})}

	\newcommand{\PSLZ}{PSL(2,\mathbb{Z})}

	\newcommand{\abs}[1]{\left\vert#1\right\vert}
	\newcommand{\norm}[1]{\left\Vert#1\right\Vert}

\begin{document}
\begin{frontmatter}

\title{Distribution of cusp sections in the  Hilbert modular orbifold}
\author{Samuel Estala Arias} 
 \ead{samuel@matcuer.unam.mx}
\address{Unidad Cuernavaca del Instituto de Matem\'aticas, Universidad Nacional Aut\'onoma de M\'exico.
\\ Av. Universidad s/n,  Colonia Lomas de Chamilpa, C.P. 62210, Cuernavaca, Morelos, M\'exico
}
 \tnotetext[t1]{Research partially supported by CONACyT grant number 129280. }
 \tnotetext[t1]{ This article constitutes the author's  Ph. D. Thesis \cite{Sam}. }

\begin{abstract}
Let $\KK$ be a number field, 
let $\mathcal{M}$ be the Hilbert modular orbifold of  $\KK$, and let  $m_q$ be the probability measure uniformly supported on  the cusp cross sections of $\mathcal{M}$  at height $q$.
We generalize a method of Zagier and show that $m_q$ distributes uniformly with respect to the normalized  Haar measure $m$ on  $\mathcal{M}$ as $q$ tends to zero, and relate  the rate by which $m_q$  approaches   $m$ to the Riemann hypothesis for the Dedekind zeta function of $\KK$.
\end{abstract} 

\begin{keyword}
Hilbert modular orbifold \sep Eisenstein series \sep Dedekind zeta function 
\end{keyword}

\end{frontmatter}

%
%
%
%
%

\section{Introduction}
\noi
Let $\HH_2=\{ x+iy \in \CC : y >0 \}$ be the Poincar\'{e} upper half-plane with the hyperbolic metric $ds^2=(dx^2+dy^2)/2$. The group $\PSL$ acts on $\HH_2$    by fractional linear transformations  which, as well, are    hyperbolic isometries.
The modular group  $\Gamma=\PSLZ$ is a discrete subgroup of  $\PSL$ and 
the quotient space $\HH_2/\Gamma$ is the classical modular orbifold.
From the classification of horocycles it follows that for each $y>0$, the modular orbifold has a unique closed horocycle $\mathcal{C}_y$  of length $y^{-1}$. 
Let  $m_y$ be the probability measure uniformly supported on $\mathcal{C}_y$ (w.r.t. arc length). Let $m$ be the normalized hyperbolic measure of $\HH_2/\Gamma$. We have the following well-known  result due to D. Zagier (cf. \cite{DZ}).
\begin{thm} \label{zagier1}Let $f$ be a smooth function on $\HH_2/\Gamma$ with compact support. Then,
$$m_y(f)=m(f)+o(y^{1/2-\epsilon}) \qquad(y \rightarrow 0)$$
for all $0<\epsilon<1/2$.
\end{thm}
\hfill

\noi In particular, the measures $m_y$  converge vaguely to  $m$  as $y  \rightarrow 0$. Besides, in \cite{DZ}   Zagier   establishes the following   remarkable equivalence to  the Riemann hypothesis. 
\begin{thm}The  Riemann hypothesis holds if and only if for every smooth function $f$ with compact support on $ \HH_2/\Gamma$ one has $$m_y(f)=m(f)+o(y^{3/4-\epsilon}) \qquad(y \rightarrow 0)$$ for all $0<\epsilon<3/4$.
\end{thm}

\hfill

\noi
In the related work \cite{Sa},  P. Sarnak proved an analogue of these theorems for the unit tangent bundle of  the modular orbifold,  which can be 
identified with the quotient space  $\PSL/\Gamma$. Likewise, in \cite{V}   A. Verjovsky  has shown that the analogue estimate of Theorem \ref{zagier1} in Sarnak's work
 is optimal  for certain characteristic functions on   $\PSL/\Gamma$. 
\hfill

\noi
The purpose of this article is to apply Zagier's theory to the case of a  number field.  Let us briefly describe  our results. Given a number filed  $\KK$ of degree $n=r_1+2r_2$ with ring of integers $\OO$, there exists
 a Riemannian manifold $\XX=(\HH_2)^{r_1} \times(\HH_{3})^{r_2} $ where the Hilbert modular group 
$\Gamma=PSL(2,\OO)$ acts  properly and discontinuously. 
If the field $\KK$ has class number $h$, then the Hilbert modular orbifold $\MM=\XX/\Gamma$ has $h$ cusps and 	
each cusp can be parametrized by the standard cusp at infinity. 
From the geometry of $\MM$, it follows that, for each cusp $\cu$ and $q>0$,
there exist a \emph{generalized closed horosphere} $B(q,\cu)/\Gamma_{\cu}$ of dimension $2r_1+3r_2-1$ and volume $q^{-1}c,$ where $c$ is  a certain constant depending on the field $\KK$. 
Let $d\upsilon_{\cu'}$ be the volume element of $B(q,\cu)$ and
$m(q,\cu)$ be the probability measure in $\MM$ uniformly supported on $B(q,\cu)/\Gamma_{\cu}$ with respect to $d\upsilon_\cu'$.
Let $m$ be the normalized Haar measure of $\MM$ and  $m_q=(m(q,\cu_1)+...+m(q,\cu_h))/h$.  We have
\begin{thm} \label{thm1} Let $f$ be a smooth function on $\MM$ with compact support. Then,
$$m_q(f)=m(f)+o(q^{1/2-\epsilon})  \qquad(q \rightarrow 0)$$
for all $0<\epsilon<1/2$.
\end{thm}
\hfill

\noi
We also give a generalization  of  Zagier's equivalence to the Riemann hypothesis and prove the following assertion.
\begin{thm}\label{thm2}
The Riemann hypothesis for the Dedekind zeta function of $\KK$ holds if and only if for every smooth 
function $f$ with compact support on $\MM$ one has
$$m_q(f)=m(f)+o(q^{3/4-\epsilon})  \qquad (q \rightarrow 0 )$$ 
for all $0<\epsilon<1/4$.
\end{thm}
\noi
We prove these theorems following the work of Zagier and Sarnak. 
We study the Mellin transform of $q^{-1}m_q$  and relate it to the Eisenstein series of the Hilbert modular orbifold by means of the Rankin-Selberg theorem. The Eisenstein series, being divisible by the Dedekind zeta function $\zeta_{\KK}(s)$ of the field,  enables us to translate  properties of $\zeta_{\KK}(s)$   into properties of the Mellin transform of $q^{-1}m_q$ and, hence, properties of $m_q$.

\section{The Hilbert Modular Group} \label{Hilbert_orbifold}
\noi In this section, we present the theory of the Hilbert modular group. 
We  quote the classical  books by Siegel \cite{S}, van der Geer \cite{GV}, Freitag \cite{F}, Elstrodt \cite{ELS}  and the article by Weng \cite{LW} for a more comprehensive introduction to Hilbert modular orbifolds.   

\hfill

\noi
Let $\HH_3=\left\{\, z=(x,y) : x \in \CC, \ y >0\, \right\}$  be the  upper half-space in Euclidean three-space with  the hyperbolic metric  $ds^{2}(z)=(d\abs{x}^2+dy^2)/y^{2}$. As usual we   identify a point
 $z\in \HH_3$ with an  element of Hamilton's  quaternions and write $
z=(x,y)=x+yj 
$. The group $PSL(2,\CC)$   acts by isometries   on $\HH_3$ as follows: for a point $z \in \HH_3$
  we have
\begin{equation} \label{accion3D}
  z \mapsto  g \cdot z =(az+b)(cz+d)^{-1} \qquad
g=\pm\begin{bmatrix}a & b\\
c&d 
\end{bmatrix}
    \notag
\end{equation}
The hyperbolic  volume element is equal to $dxdy/y^3$ where $dx$  is the two dimensional Lebesgue measure on $\CC$.

\hfill

\noi
Let $\KK$ be an algebraic number field  with $r_1$ real places and  $r_2$ complex  places and write $n=[K:\QQ]=r_1+2r_2$ for the degree of $\KK$ over $\QQ$. Consider the space
$\XX=\HH_2^{r_1} \times\HH_3^{r_2}$, and denote a point in $\XX$ by $\zz=(z_1,..., z_r)$, where $r=r_1+r_2$. The space $\XX$ has   the  Riemannian metric $$
ds^2(\zz)=\sum_{i=1}^{r}( d\abs{x_i}^2+dy_i^2) /y_i^{2}
$$ with corresponding  volume element
$
d\upsilon(\zz)=d\xx d\yy/y_1^2\cdots y_{r}^3,
$ where $d\xx d\yy$ is the Lebesgue measure on $(\RR^{r_1}\times\CC^{r_2})\times\RR^{r}$.
\noi  
The group $G=PSL(2,\RR)^{r_1} \times PSL(2,\CC)^{r_2}$ acts transitively on $\XX$ by isometries:
$$g \cdot \zz:=\big(g_1\cdot z_1,...,g_r \cdot z_r\big)$$ 
where $g=(g_1,...,g_r) \in G $.
The isotropy group at the point $(i,...,j)\in \XX$ is the (maximal) compact subgroup $PSO(2,\RR)^{r_1}\times PSU(2)^{r_2}$ and we have a natural identification between $\XX$ and the quotient of $G$ by this subgroup.

\hfill

\noi
To the  real and complex places of $\KK$ there correspond $r$ embeddings of $\KK$ into  $\CC$. We arrange these embeddings as usual, and write them as  $\alpha^{(i)}$, for $\alpha \in \KK$. Alike, we often write $\alpha=\alpha^{(1)}$.
The group $PSL(2,\KK)$ can be embedded in $G$ as follows
$$
\pm \begin{bmatrix}
a & b \\
c & d 
\end{bmatrix} \mapsto  \left(  \pm\begin{bmatrix} a^{(1)} & b^{(1)} \\
c^{(1)} & d^{(1)}
\end{bmatrix}
,...,\pm \begin{bmatrix} a^{(r)} & b^{(r)} \\
c^{(r)} & d^{(r)}
\end{bmatrix}  \right).
$$ 
Let $\OO$ be the ring of algebraic integers of $\KK$. The \emph{Hilbert modular group} for the field $\KK$ is the group $\Gamma=PSL(2,\OO)$ contained in $PSL(2,\KK) \subset G$. The action of $\Gamma$ in $\XX$ is irreducible and,  since $\OO$ is a lattice, $\Gamma$ is a discrete subgroup of $G$.
 Then a theorem of Thurston \cite[Ch. 13]{Th} implies that $\Gamma$  acts properly and discontinuously on $\XX$ and that
the quotient space $$\MM=\Gamma\backslash \XX$$ is 
a differentiable orbifold. It is  called the \emph{Hilbert modular orbifold} of the field $\KK$.

\hfill

\noi
The action of $G$ in $\XX$ extents  to the ideal boundary $\mathbb{P}=\mathbb{P}(\RR)^{r_1} \times\mathbb{P}(\CC)^{r_2}$ by fractional linear transformations. The cusps of $\Gamma$ are  the parabolic fixed points of $\Gamma$ in $\mathbb{P}$. They are the image of the application $\mathbb{P}(\KK) \rightarrow\mathbb{P}$ given by $\alpha/\beta \mapsto (\alpha^{(i)}/{\beta^{(i)}})$. The group $\Gamma$ acts on the set of cusps, and a classical theorem of Maass, shows that non equivalent cusps are in one-to-one correspondence with the elements of the ideal class group of $\KK$. To the class of a fractional ideal $\fa=\la\sigma, \rho\ra \subset \KK$ there corresponds
the class of the cusp $\lambda=\sigma/\rho$. In particular, to the identity in the ideal class group of $\KK$, that is to say, to the class of principal ideals, there corresponds the class of the cusp at infinity $(\infty,...,\infty)\in \mathbb{P}$. Then, if   $h$  is the class number of  $\KK$, we know that $\Gamma$ has $h$ inequivalent cusps.

\hfill

\noi
Let $\lambda\in \mathbb{P}(\KK)$ be a cusp  and write $\lambda=\sigma/\rho$, $\fa=\langle \sigma,\rho\rangle$. Since $\rho$ and $\sigma$ generate 
$\mathfrak{a}$, there exist $\xi, \eta \in\mathfrak{a}^{-1}$ such that $\rho \eta-\sigma \xi=1$,  and  we have a matrix
$$
\mathrm{A}=\begin{bmatrix}
\rho & \xi \\
\sigma & \eta \end{bmatrix}
 \in SL(2, \KK)
$$ with the property $\mathrm{A}^{-1}(\lambda)=\infty$. The matrix $\mathrm{A}$ is called a \emph{matrix associated to} $\lambda$. 
 This matrix depends on the representation of $\lambda$ and the choice of $\xi,\eta$. However, any other matrix associated to $\lambda$ as above  is  equal to $$
\mathrm{A}\begin{bmatrix}
\alpha & \beta \\
0 & \alpha^{-1} \end{bmatrix}$$ for some $\alpha, \beta \in\KK$ . 
\begin{prop}Let $\lambda \in \mathbb{P}(\KK)$ be a cusp and $\Gamma_{\lambda}=\left\{\, \gamma \in \Gamma  \mid  \gamma\lambda=\lambda \, \right\}$ be  the isotropy group of $\Gamma$ at $\lambda$. Then, the group $\Gamma_\lambda$ is described as follows
\begin{equation}
\Gamma_{\lambda}= \left\{ \,\mathrm{A} \begin{bmatrix}\varepsilon & \zeta\varepsilon^{-1}\\ 0 & \varepsilon^{-1} \end{bmatrix}
\mathrm{A}^{-1}  \mid \varepsilon \in \OO^{\times}, \quad \zeta \in \mathfrak{a} ^{-2}\, \right\} \notag
\end{equation}
where $\OO^{\times}$ denotes the invertible elements of $\OO$.
\end{prop}
\hfill

\noi
Let us describe the action of $\Gamma_{\lambda}$ on $\XX$. First, the fractional ideal
$\mathfrak{a}^{-2}$ is a free $\ZZ$-module of rank $n$ and, therefore, has $\ZZ$-bases.  Let
 $\alpha_1,...,\alpha_n$ be such a base and let $\NN(\fa)$ be the norm of the ideal $\fa$.  Likewise, from the Dirichlet's  units theorem, there exist \emph{fundamental units} $\varepsilon_{1},..., \varepsilon_{r-1} \in \OO^{\times}$, such that 
$$
\OO^{\times}=\{  \varepsilon_1^{k_1} \cdots \varepsilon_{r-1}^{k_{r-1}} \ | \ k_1,...,k_{r-1} \in \ZZ \ \}\times W,
$$
where $W$ is the group of roots of unity contained in $\KK$.
Let $\zz=(z_1,...,z_r)$ be any point in $\XX$ and define $\zz^{\star}=\mathrm{A}^{-1}\zz$. We write
$$
\zz^{\star}=(z_1^{\star},..., z_r^{\star}), \quad \xx^{\star}=(x_1^{\star},...,\Re(x_r^{\star}),\Im(x_r^{\star})), \quad \text{ and } \yy^{\star}=(y_1^{\star},...,y_r^{\star}).
$$ Set  $\NN(\yy^{\star})=y_1^{N_1} \cdots y_r^{N_r}$, where $N_i$ is the degree of the embedding $i$.
The local coordinates of $\zz$ at $\lambda$ are defined to be the $2r$ quantities
\begin{equation}
q, \Y_1,...,\Y_{r-1}, \X_1,...,\X_{r} \notag
\end{equation}
determined by the relationships (compare to \cite{LW}):
\begin{align}
q&=\NN(\fa)\NN(\yy^{\star})\notag \\
\begin{bmatrix}
\ln \abs{\varepsilon_1^{(1)}} & \cdots &\ln \abs{\varepsilon_{r-1}^{(1)}} \\
\vdots & \ddots &\vdots\\
\ln \abs{\varepsilon_1^{(r-1)}} & \cdots &\ln \abs{\varepsilon_{r-1}^{(r-1)}}
\end{bmatrix}\begin{bmatrix}\Y_1\\ \vdots\\ \Y_{r-1}
\end{bmatrix}&=\begin{bmatrix}\frac{1}{2}\ln\frac{y_1^{\star}}{\sqrt[\leftroot{2} \uproot{2} n]{\NN(\yy^{\star })}}\\ \vdots \\\frac{1}{2}\ln\frac{y_{r-1}^{\star}}{\sqrt[\leftroot{2} \uproot{2} n]{\NN(\yy^{\star })}} 
\end{bmatrix}  \notag \\
 \qquad \begin{bmatrix}
                   \alpha_{1}^{(1)}& \cdots & \alpha_{n}^{(1)} \\
\vdots & \ddots & \vdots \\
\Re(\alpha_{1}^{(r)}) & \cdots & \Re(\alpha_{n}^{(r)}) \\
\Im(\alpha_{1}^{(r)}) & \cdots & \Im(\alpha_{n}^{(r)})
                  \end{bmatrix} \begin{bmatrix} \X_1 \\ \vdots \\ \Re(\X_r)\\ \Im(\X_r)
\end{bmatrix} &= \begin{bmatrix} x_1^{\star} \\ \vdots \\ \Re(x_r^{\star})\\ \Im(x_r^{\star})
\end{bmatrix}\notag
\end{align}

\begin{lemma}\label{lemadelhaz}
\noindent Let $\zz\in \XX$ be a point, $\zeta = m_1\alpha_1+\cdots +m_n\alpha_n\in \mathfrak{a}^{-2}$ and $\varepsilon= w \varepsilon_{1}^{k_1} \cdots \varepsilon_{r-1}^{k_{r-1}} \in \OO^{\times}$, with $w\in W$. Then, under the modular transformation $$\mathrm{M}=\mathrm{A}
\left[\begin{matrix}\varepsilon & \varepsilon^{-1}\zeta \\ 0 & \varepsilon^{-1}
\end{matrix}
\right] \mathrm{A}^{-1}$$ the local coordinates of $\zz$ become
$$
q, \Y_1+k_1,...,\Y_{r-1}+k_{r-1}, \X_1^{\star}+m_1,...,\Re(\X_r^{\star})+m_{n-1}, \Im(\X_r^{\star})+m_n
$$
where the column vector $\X^{\star}=(\X_1^{\star},...,\Re(\X_r^{\star}),\Re(\X_r^{\star}))$
is described as follows: if we write the definition of the local coordinates in matrix notation as $\mathrm{O}(\X)=\xx^{\star}$ with the column vectors $\X=(\X_1,..., \Re(\X_r), \Im(\X_r))$  and  $\xx^{\star}$, then $\X^{\star}=\mathrm{O}^{-1}\mathrm{E}^2\mathrm{O} \X $ where $\mathrm{E}$ is the block matrix  $$
\mathrm{E}=\left[\begin{matrix}\varepsilon^{(1)} & ...& 0& 0 \\
\vdots & \ddots & 0 & 0 \\ 
0 & \cdots&\Re(\varepsilon^{(r)}) &  -\Im(\varepsilon^{(r)})\\
0 & \cdots&\Im(\varepsilon^{(r)}) &\Re(\varepsilon^{(r)})
\end{matrix}
\right]
$$ 

\end{lemma} 

\begin{pf}
This follows  from 
\begin{align}(\mathrm{M}\zz)^{\star}=\mathrm{A}^{-1}\mathrm{MA}\zz^{\star}=(\varepsilon^{2} \xx^{\star}+\zeta,\abs{\varepsilon}^{2}\yy^{\star})=\big((\varepsilon^{(i)})^{2} x_i^{\star}+\zeta^{(i)},\abs{\varepsilon^{(i)}}^{2} y_i^{\star}\big). \notag
\end{align} 
\end{pf}
\hfill
\begin{remark}
If  the group $W$ of roots of unity is different from $\pm 1$, then all Galois embeddings of $\KK$ are complex. Let $w \in W$ be different from $\pm 1 $, then the action of $\mathrm{A}
\left[\begin{matrix} w &  0\\ 0 & w^{-1}
\end{matrix}
\right] \mathrm{A}^{-1}$ in $\XX=(\CC\times \RR_{+})^{r_2}$ 
can be written in the coordinates $(\xx^{\star},\yy^{\star})=\mathrm{A}^{-1}\zz$ as  the  transformations $(x_i^{\star},y_i^{\star}) \mapsto \big((w^{(i)})^{2}x_i^{\star},y_i^{\star}\big)$, where each $w^{(i)}$ is a root of unity.  
\end{remark}

\hfill
\begin{definition}
A point $\zz=(\xx,\yy)$ in $\XX$ is called $\mathsf{reduced}$ $\mathsf{with}$ $\mathsf{respect}$ $\mathsf{to}$ $\lambda$  if its local coordinates at $\lambda$ satisfy
$$0<q, \quad
-\frac{1}{2} \leq \Y_1,...,\Y_{r-1},\X_1,...,\Re(\X_{r}),\Im(\X_r) < \frac{1}{2}$$
and his coordinate $\xx^{\star}$ belongs to a fundamental domain of the action of $W$ in $\CC^{r_2}$ by $\xx^{\star} \mapsto w^{2}\xx^{\star} $.
\end{definition}
\hfill

\noi
As a result  of lemma \ref{lemadelhaz} the set $\mathscr{F}_{\lambda}$ of reduced points (w.r.t. $\lambda$ ) in $\XX$ is a fundamental domain for $\Gamma_{\lambda}$. Moreover, define 
the \emph{generalized horosphere} at height $q$ (or distance to infinity $1/q$) of the cusp $\lambda$  by $B(q,\lambda)=\{ \zz \in \XX \, \mid \, \mu(\lambda, \zz) = q\}$. Then,  
the action of $\Gamma_{\lambda}$  on $\HH$ reduces to its action on $B(q,\lambda)=\{ \zz\in \XX \,\mid \, \mu(\zz,\lambda)=q\}$. 

\hfill

\begin{lemma} Let  $d\upsilon'_\lambda$ be  the measure induced on $B(q,\lambda)$ by the Riemannian metric on $\XX$ and let $d\upsilon$ be the Haar measure on $\XX$. 
Let $\R$ be the regulator  and let $\D$ be the absolute value of the discriminant of $\KK$.
Then, in the local coordinates of $\lambda$, we have
$$d\upsilon'_{\lambda}(\zz)=(r_1+4r_2)^{1/2}\,2^{r_1-r_2-1}\,q^{-1}\,\sqrt{\D}\, \R \,d\X d\Y $$ 
where $d\X d\Y$ is the Lebesgue measure in $\RR^{r_1}\times \CC^{r_2} \times \RR^{r-1}$. Likewise, we have 
$$d\upsilon(\zz)= 2^{r_1-r_2-1}\,q^{-2}\,\sqrt{\D}\, \R \,d\X d\Y dq $$ 
\end{lemma}
\begin{pf}
 Since $\mathrm{A}$ is an isometry,  Riemannian  metric objects   in the coordinates $\zz$ can be described in the coordinates $\zz^{\star}$. Now, since the change between the $\X$'s and the $\xx^{\star}$'s is linear,  we have $d\X_1\cdots d\X_r=(2^{r_2}\NN(\fa)^{2}/\sqrt{\D})\, d\xx^{\star}$ (see \cite{JN}). Then, the  volume induced by the Riemannian metric on the submanifold $(\xx^{\star},\yy^{\star})$, with a fixed $\yy^{\star}$, is given by
$$
\frac{d\xx^{\star}}{\NN(\yy^{\star})}=q^{-1}\,2^{-r_2}\, \sqrt{\D}\ d\X_1\cdots d\X_d
$$ On the other hand, from the definition of the local coordinates, if $\hat{\q}=\NN(\yy^{\star})$, we have
\begin{align}
 y_i^{\star}&= \hat{q}^{\frac{1}{n}} \exp\big(\sum_{k=1}^{r-1} 2\Y_k \ln \abs{\varepsilon_k^{(i)}}\big)
=\hat{q}^{1/n}\,1\prod_{k=1}^{r-1}\abs{\varepsilon_{k}^{(i)}}^{2\Y_k} \notag
\end{align} for all $i=1,...,r$. 
The application $(y_i^{\star})\mapsto (\log y_i^{\star})$ is an isometry between  $(\RR_{+}^{r},\prod_{i=1}^{r}dy_i^{\star}/y_i^{\star})$ and $(\RR^{r},\sum_{i=1}^{r}dt_i)$. The local coordinates of $\yy^{\star}$ factor
through a commutative diagram
\begin{equation}\begin{CD}
\RR^{r-1}\times\RR_{+} @>>>\RR^{r}_{+} \\
@V{(\Y,q) \mapsto (\Y,\log q)}VV @VV{\log}V\\
\RR^{r} @>f>> \RR^{r} \notag
\end{CD}\end{equation}
where $\big(f(\Y_1,...,\Y_{r-1},\log q)\big)_i=(\sum_{k=1}^{r-1}2\Y_k\log\abs{\varepsilon_k^{(i)}} +\frac{1}{n}\log q)$ is a linear transformation. 
Then, the measure induced on the submanifold $\NN(\yy^{\star})=\hat{q}$ by the measure $\prod_{i=1}^{r}dy_i^{\star}/y_i^{\star}$   is equal to ( cf. \cite{JN}) 
$$
(r_1+4r_2)^{1/2}\,2^{r_1-1}\,\R \,d\Y_1\cdots d\Y_{r-1},
$$
where $d\Y_1\cdots d\Y_{r-1}$ is the Lebesgue measure on $\RR^{r-1}$. Moreover,
$$
q^{-1}\,2^{r_1-1}\,\R \,d\Y_1\cdots d\Y_{r-1}dq=\prod_{i=1}^{r}dy_i/y_i.
$$
From these results the lemma follows.
\end{pf}
\hfill

\begin{remark}For each $q>0$ the action of $\Gamma_{\lambda}$ leaves $d\upsilon'_{\lambda}$ invariant and the quotient $\Gamma_{\lambda}\backslash B(q,\lambda)$ has volume $$q^{-1}\,\omega^{-1}\,2^{r_1-r_2}\,(r_1+4r_2)^{1/2} \,\sqrt{\D}\,\R,$$ where $\omega$ is the number of roots of unity contained in $\KK$. Alike,
we have a covering orbifold application $i_\lambda:\Gamma_{\lambda}\backslash \XX \rightarrow \MM$, which gives the immersion $i_{\lambda}\big(\Gamma_{\lambda}\backslash B(q,\lambda)\big) \rightarrow \MM$. 
\end{remark}
\hfill

\noi
We now construct a fundamental domain for the action of $\Gamma$ in $\XX$ and describe the topology of the quotient. Let $\lambda=\frac{ \rho}{\sigma} \in \mathbb{P}(\KK)$ be a cusp and write $\mathfrak{a}=(\rho,\sigma)$. For $\zz \in \XX$ we have defined the height of $\zz$ at $\lambda$ as
\begin{align}
\mu(\lambda,\zz) &=\frac{\NN(\mathfrak{a})^2 \,\NN(\yy)}{\abs{\NN(-\sigma \zz+\rho)}^2} \notag \\
\notag
\end{align} This is independent of the choice of $\rho$ and $\sigma \in \KK$ and, for any $\gamma \in \Gamma, $  we have the invariance property $
\mu(\gamma(\lambda),\gamma(\zz))=\mu(\lambda,\zz)
$.  
The proofs of the following propositions can be found in  \cite{LW}.
\begin{prop} \label{independencia_cuspidal}
There exists a positive number $l_1,$ depending only on $\KK,$ such that for $\zz \in  \XX$ the inequalities $\mu(\lambda,\zz)>l_1$ and $\mu(\tau,\zz)>l_1$ for  $\lambda,\tau \in \mathbb{P}(\KK)$
imply $\lambda=\tau$.
\end{prop}

\begin{prop}\label{cota_universal} There exists a positive  number $l_2,
$ depending only on $\KK,$ such that for $\zz \in \XX$ there exists a cusp $\lambda$ such that $\mu(\lambda,\zz)>l_2$.
\end{prop}
\hfill

\noi
Let $\lambda$ be a cusp. The \emph{sphere of influence} of $\lambda$ is defined by$$
S_{\lambda}=\left\{\, \zz \in \XX \, \mid \, \mu( \lambda, \zz) \geq  \mu (\tau,\zz ) \ \ \forall \tau \in \mathbb{P}(\KK) \,\right\}
$$
The invariance condition $\mu(\gamma \zz,\gamma \lambda)=\mu(\zz, \lambda)$ implies $S_{\gamma (\lambda)}=\gamma(S_\lambda)$. The boundary of $S_{\lambda}$ consists of pieces  defined by equalities $\mu(\lambda, \zz)= \mu( \tau, \zz)$ with  $\lambda \neq \tau$. Moreover, the action of $\Gamma$ in the interior $F_{\lambda}^{0}$ of $S_{\lambda}$ reduces to that of the isotropy group
$\Gamma_{\lambda}$ at $\lambda$, i.e, if $\zz$ and $\gamma \zz$ both belong to $F_{\lambda}^{0}$, then $\gamma \lambda= \lambda$. 
\noi
Let $i_{\lambda}: \Gamma_{\lambda} \backslash  S_{\lambda} \rightarrow \Gamma \backslash \XX$ be the natural map. From   propositions \ref{independencia_cuspidal}  and \ref{cota_universal}, we have
$$
\Gamma \backslash \XX = \cup_{\lambda}i_{\lambda}(\Gamma_{\lambda} \backslash S_{\lambda} )
$$
where the union is taken over a set of $h$ non equivalent cusps.

\hfill

\noi
Let $G_{\lambda}$ be the intersection of the sphere of influence at $\lambda$ and the reduced points with respect to $\lambda$. Then $G_{\lambda}$ is a fundamental region for the action of $\Gamma_{\lambda}$ in $S_{\lambda}$.  
Since $\MM$ is connected, there exist
$\lambda_1,...\lambda_{h}$ non equivalent cusps such that $\cup_{i=1}^{h} G_{\lambda_i}$ is connected and  a fundamental domain for  $\Gamma$.
For $\T$ sufficiently large   we  have a compact orbifold with boundary:
$$
\MM_{\T}=\Gamma \backslash \left\{ \zz \in  \XX: \mu(\zz,\lambda) \leq  \T \text{ for all cusps } \lambda \right\}
$$  
The boundary of $\MM_{\T}$ consists of $h$  compact orbifolds $i_{\lambda} ( \Gamma_{\lambda} \backslash B(\lambda, r ))$ of dimension $2r_1+3r_2-1$ and, as a topological space, $\Gamma\backslash \XX$ has $h$ ends: $$ \Gamma\backslash \XX= \MM_{\T} \cup_{\substack{\partial \MM_{\T}}} \big(\partial \MM_\T \times [0, \infty)\big).$$

%
%

%
%

\section{Non holomorphic  Eisenstein series.}
\noi
We now define the Eisenstein series  of the Hilbert modular orbifold and state their properties. Investigations on the Eisenstein series have been quite extensive. We quote
the  classical references  \cite{Ku}, \cite{He},\cite{Ef} and the  exposition \cite{CP}.

\hfill

\noi
First, recall that $\XX$ (as a Riemannian manifold) have associated the Laplace-Beltrami operator:
\begin{equation}\label{laplaciano}
\Delta=\sum_{i=1}^{r_1}y_i^2\big(\frac{\partial^{2}}{\partial x_i ^2}+\frac{\partial^{2}}{\partial y_i^2}\big)
+\sum_{i=r_1+1}^{r}y_i^2\big(\frac{\partial^{2}}{\partial x_i \partial \overline{x_i}}+\frac{\partial^{2}}{\partial y_i^2}\big)-y_i\frac{\partial}{\partial y_i} \notag
\end{equation}
It is invariant under the action of $G$ on functions; that is, for all  $g\in G$, we have
$
\Delta(f\circ g)=\Delta (f) \circ g,
$ for any smooth function $f$ on  $\XX$.  
In particular, since $$
 \Delta(\NN(\yy)^s)=(r_1+4r_2)s(s-1)\NN(\yy)^s
$$ we have $$ \Delta(\mu(\lambda, \zz) )= (r_1+4r_2)s(s-1)\mu(\lambda, \zz)$$ for all points $\lambda\in \mathbb{P}(\KK)$.

\hfill

\noindent
Let $\cu$ be a cusp of the Hilbert modular group. The Hilbert modular Eisenstein series associated to $\cu$ is defined by the series 
\begin{equation}\label{Emodular} \E_{\cu}(\zz,s)=\sum_{\gamma \in \Gamma_{\cu} \backslash \Gamma} \mu(\cu,\gamma \zz)^s \notag
\end{equation}
where  the sum
is taken over any complete set of  representatives of the  classes $\Gamma_{\cu} \backslash \Gamma$. This is independent of  $\cu$ in an equivalence class of cusps. The series  defining $\E_{\cu}(\zz,s)$ converges absolutely and uniformly for $s$ in bands $\{\sigma_1>\Re(s)>\sigma_{0}>1 \}$ and  $\zz$ in every compact subset of $\XX$. Therefore,  it defines a continuous function which is holomorphic  on  the  half-plane $\Re(s)>1$. Moreover, it represents  a $\Gamma$ automorphic form, that is
$$
\E_{\cu}(\zz,s)=\E_{\cu}(\gamma \zz,s) \quad \text{ for all } \gamma \in \Gamma
$$ and satisfies   the differential equation
\begin{equation}\label{ecuacion_armonica}
\Delta \E_{\cu}(\zz,s)=(r_1+4r_2)s(s-1)\E_\cu( \zz,s) \notag
\end{equation}

\noindent
Recall that $\OO^{\times}$ acts on the set of pairs $(\alpha,\beta)\in \KK^{2}$ by $\varepsilon\cdot(\alpha,\beta)=(\varepsilon\alpha,\varepsilon\beta)$ and  that two pairs in the same orbit are called associated. Let  $\{\mathrm{M}_j\}_{j\in \NN}$ be a complete set of representatives for the quotient $\Gamma_{\lambda} \backslash \Gamma$ and, let $\mathrm{A}$ be a matrix associated to $\lambda$. Write
$$\mathrm{A}^{-1}\mathrm{M}_{j}=\left[\begin{matrix}\alpha_j & \beta_j \\ \gamma_j & \delta_j
\end{matrix}
\right] \qquad \forall j \in \NN 
$$ Then,  the collection  $\{(\gamma_j,\delta_j): j \in \NN\}$ is a complete set of non associated pairs of generators for the ideal $\mathfrak{a}$  and, from the definition of $\mu(\lambda,\zz)$, we have
$$
\E_{\cu}(\zz,s)=\NN(\mathfrak{a})^{s}\sum _{\substack{(\gamma, \delta) \in \mathfrak{a}^2 /\OO^{\times} \\ (\gamma, \delta)=\mathfrak{a}}} \frac{\NN(\yy)^{s}}{\abs{\NN(\gamma \zz+\delta)}^{2s}}.
$$
In order to complete the Eisenstein series, we need the partial zeta function of an ideal class $\mathcal{A}$. It is
defined by the series
$$
\zeta_{\KK}(s ,\mathcal{A})= \sum_{\substack{ \mathfrak{a} \in \mathcal{A} \\  \text{ integral }}} \frac{1}{\NN(\mathfrak{a})^s}.
$$
The series is absolutely convergent for $\Re(s)> 1$, and uniformly in $\Re(s)>1+\epsilon$. We define the completed zeta function of an ideal class $\mathcal{A}$ by $\Lambda(s)\zeta_{\KK}(s,\mathcal{A})$ where, as in \cite{SL},
$$\Lambda(s)=2^{-r_2 s} \D^{\frac{s}{2}} 
\pi^{-\frac{ns}{2}}\Gamma\left(\frac{s}{2} \right)^{r_1} \Gamma(s)^{r_2}.
$$  
By a theorem of Hecke, $\zeta_{\KK}^{\star}(s,\mathcal{A})$ has a meromorphic continuation to the whole complex plane and satisfies the functional equation
$$
\zeta_{\KK}^{\star}(s,\mathcal{A})=\zeta_{\KK}^{\star}(1-s,\mathcal{A}'),
$$
where $\mathcal{A}\mathcal{A}'= [\mathfrak{d}]$ is the ideal class of the different ideal of the field $\KK$. Recall that, by definition, 
$\mathfrak{d}^{-1}$ is the dual module of $\OO$, w.r.t. the trace function of $\KK$. 

\hfill

\noindent
For a cusp $\tau$ denote by $\mathcal{A}_\lambda$ the ideal class associated with $\tau$. Let $\lambda$ be a cusp,  the associated normalized Hilbert modular Eisenstein series  is
$$
\E_{\cu}^{\star}(\zz,s)= \sum_{\substack{\cu' }} \zeta_{\KK}^{\star}(2s,\mathcal{A}^{-1}_{\cu}\mathcal{A}_{\cu'})
\E_{\cu'}(\zz,s) \qquad (\zz \in  \XX, \ \Re(s)>1)
$$ where the sum is taken over a set of non equivalent  cusps1.
Writing the integral ideals of an ideal class $\mathcal{C}=[\mathfrak{c}]$ as the set of $\xi \mathfrak{c}$ with $\xi \in \mathfrak{c}^{-1}/\OO^{\times}$, one can show that
$$
\E_{\cu}^{\star}(\zz,s)=\Lambda(2s)\NN(\mathfrak{a})^{s}\sum _{\substack{(\gamma, \delta) \in \mathfrak{a}^2 /\OO^{\times} \\ (\gamma, \delta)\neq 0}}\frac{\NN(\yy)^{s}}{\abs{ \NN(\gamma \zz+\delta)}^{2s}}.
$$

\noindent
Let $\cu'\in \mathbb{P}(\KK)$ be another cusp and write $\cu'=\sigma'/\rho'$,  $\fa'=\la\sigma' \rho'\ra$. Let
$$\mathrm{A}'=\left[\begin{matrix}\rho' & \xi' \\ \sigma' & \eta'
\end{matrix}
\right] $$ be a matrix associated to $\cu'$. Since  Eisenstein series are  $\Gamma$-automorphic functions, in particular, we have that  
$$
\E_{\cu}^{\star}(\mathrm{A}'\zz,s)=\Lambda(2s)\NN(\mathfrak{a})^{s}\sum _{\substack{(\gamma, \delta) \in \KK^2 /\OO^{\times} \\0 \neq (\gamma\eta'-\delta\sigma',\gamma\xi' -\delta \rho')\subset \fa}}\frac{\NN(\yy)^{s}}{\abs{ \NN(\gamma \zz+\delta)}^{2s}}
$$
is invariant for the lattice $\fb'=\fa'^{-2}\subset \RR^{r_1}\times\CC^{r_2}$. Hence, it has  Fourier expansion
$$
\E_{\cu}(\mathrm{A}'\zz,s)=\sum_{\substack{l\in \fb'^{\star}}}a_l^{\cu,\cu'}(\yy,s)e^{2\pi i Tr (l\xx)}
$$ 
where $\mathfrak{b}'^{\star}=\left\{ \xi \in \KK \, | \, Tr(\xi \OO) \subset \ZZ\right\}
$ is the dual module of $\fb'$ and
$$
Tr(x)=\sum_{\substack{i \leq r_1}} x_i+ \sum_{\substack{i>r_1}}2 \Re(x_i),
$$
for  $x=(x_1,...,x_r)\in \RR^{r_1}\times\CC^{r_2}$. 
In order to express  the Fourier coefficients of the Eisenstein series we  need the MacDonald Bessel function:

$$K_s(y)=\frac{1}{2}\int_{0}^{\infty} e^{-y(t+t^{-1})/2}t^s\frac{dt}{t} \quad \text{ for }s \in \CC, \ y>0. 
$$
The Bessel function has many properties (cf. \cite{DB}): 
\begin{enumerate}[i)]
 \item  it is invariant under the transformation $s \mapsto -s$, that is $K_s(y)=K_{-s}(y)$.
 \item  all derivatives of $K_s(y)$ with respect to $y$ are of rapid decay when  $y \rightarrow+\infty$,  i.e.,
$$
\abs{K_s^{(l)}(y))} \leq K_{\Re(s)}^{(l)}(2)e^{-\frac{y}{2}} \quad \forall y >2 \qquad  l=0,1,2,...
$$ 
\item it satisfies the differential equation 
$$
\left\{ y^{2}\frac{d^2}{dy^2}+y \frac{d}{dy}-(y^{2}+s^{2}) \right\} K_{s}(y)=0.
$$
\item it is related to a Fourier transform, i.e.,  for $l\in \RR/\{0\}$ and  $\Re(s)>1$, we have
$$
y^{s}\pi^{-s}\Gamma(s)\int_{\RR}\frac{e^{2\pi i lt}}{(t^{2}+y^{2})^{s}}dt=2 \abs{l}^{s-\frac{1}{2}} \sqrt{y} K_{s-\frac{1}{2}} (2\pi  \abs{l} y) 
$$ 
\end{enumerate}
\hfill

\noi
The following theorem on  the Fourier expansion of the Eisenstein series is proved by Sorensen in \cite{CM} (see also \cite{L2}, \cite{LW}).
\begin{prop}
Let $\fb'^{\star}$ be the dual module of $\fb'=\fa'^{-2}$.  Write $\zz^{\star}=(\xx^{\star},\yy^{\star})=\mathrm{A}'^{-1}\zz$ and $q_{\cu'}=\mu(\cu',\zz)$. Then,  $\E_{\cu}^{\star}(\zz,s)$ 
has  the Fourier expansion at the cusp $\cu'$
\begin{align}
\E^{\star}_{\cu}(\zz,s)&=\zeta_{\KK}^{\star}(2s,\mathcal{A}_{\cu}^{-1}\mathcal{A}_{\cu'})q_{\cu'}^{s}+\zeta_{\KK}^{\star}(2s-1,\mathcal{A}_{\cu}^{-1}\mathcal{A}_{\cu'}^{-1})q_{\cu'}^{1-s} \notag \\
&+2^{r}q_{\cu'}^{1/2}\sum_{\substack{l\in \fb'^{\star}\\ l \neq 0}}\tau^{\cu,\cu'}_{1-2s}(l)K_{s-1/2}(\yy^{\star},l)e^{2\pi i Tr (\xx^{\star}l)} \notag \end{align} 
where
\begin{align}
\tau_{s}^{\cu,\cu'}(l)&=\NN(\fb'\mathfrak{d}l)^{-s/2}\sum_{\substack{\mathfrak{q} \in \mathcal{A}_{\cu}^{-1} \mathcal{A}_{\cu'}^{-1}\\ \text{integral} \\
\mathfrak{q}: \fb'\mathfrak{d}l }} \NN(\mathfrak{q})^{s} \notag
\end{align}
and $$K_s(\yy^{\star},l)=K_s(2\pi y_1 ^{\star}\abs{l^{(1)}})\cdots K_s(4\pi y_r^{\star}\abs{l^{(r)}})$$
\end{prop}
\hfill

\noi
Summing the partial Eisenstein series over a set of non equivalent   cusps we have
$$
\E(\zz,s)=\sum_{\substack{\cu  }}\E_{\cu}(\zz,s), \qquad \E^{\star}(\zz,s)=\sum_{\substack{\cu}}\E_{\cu}^{\star}(\zz,s)
$$ The series $ \E$ and $\E^{\star}$ are related by $
\E^{\star}(\zz,s)=\zeta^{\star}_{\KK}(2s)\E(\zz,s)$. 
By exploring the Fourier expansion of $\E^{\star}(\zz,s)$ at  the cusp $\cu'=\infty$
we can see that it represents a meromorphic function which is holomorphic in $\CC$
except for two simple poles at $s=0,1$ and satisfies the functional equation
$$
 \E^{\star}(\zz,s)= \E^{\star}(\zz,1-s) \qquad\forall \zz \in  \XX \quad \forall s \in \CC \backslash \{0,1\}.
$$ From the class number formula, the residue of $\E^{\star}(\zz,s)$ at the simple pole $s=1$ is equal to $2^{r_1-1}\,\R\, h \, \omega^{-1}$.  
Besides, counting the poles and zeros of 
$\E^{\star}(\zz,s)$, we can see that the poles of $ \E(\zz,s)$  are in the zeros of the function $\zeta_{\KK}^{\star}(2s)$. Since,  from the Euler product representation of the  Dedekind zeta function and the Landau prime number theorem, $\zeta_{\KK}^{\star}(s)$ does not vanish in the region $\Re(s)\geq 1$,  we have that $\E(\zz,s)$ is holomorphic for $\Re(s)\geq \frac{1}{2}$  except for a simple pole at $s=1$ with residue:
\begin{equation}
\mathcal{R}es_{s=1}(\E(\zz,s))=\frac{2^{r-1}\,h\,\R}{\omega \, \pi^{-n}\, \D \,\zeta_{\KK}(2)} \qquad (\forall \zz \in \XX)
\notag
\end{equation}

\begin {lemma}\label{sii} The  Riemann  hypothesis for the Dedekind zeta function $\zeta_{\KK}(s)$ holds if and only if    for all $ \zz \in \XX$ the function $\E(\zz,s)$ is holomorphic in the half plane  $\Re(s)>1/4$, except for a simple pole at $s=1$ .
\end{lemma}

\begin{pf} Suppose that 
$\E(\zz,s)$ is holomorphic in the half plane $\Re(s)>\frac{1}{4}$ for all $\zz \in \XX$ except for a simple pole at $s=1$. Let $s=\sigma+it$ be   such that $1/4\leq\sigma<1/2$ and $\zeta_{\KK}^{\star}(2s)=0$,  then   $\zeta_{\KK}^{\star}(2s-1)\neq 0$  and the zero Fourier coefficient of $\E^{\star}(\zz,s)$ does not vanish. Therefore, there exists $\zz \in \XX$ such that $\E^{\star}(\zz,s)=\zeta_{\KK}^{\star}(2s)\E(\zz,s) \neq 0$ and so $\sigma=1/4$. The converse claim is clear.
\end{pf}
\hfill

\begin{coro}\label{expansion_final}
Let $\cu$ be a cusp, let $\fa$ be the ideal associated to $\cu$ and $\fb=\fa^{-2}$. Let $q_\cu=\mu(\zz,\cu)$ and $\zz^{\star}=\mathrm{A}^{-1}(\zz)$. Then, the function $\E(\zz,s)$   has  the Fourier expansion at  the cusp $\cu$
\begin{align}
\E(\zz,s)&=q_\cu^{s}+\frac{\zeta_{\KK}^{\star}(2s-1)}{\zeta_{\KK}^{\star}(2s)}q_\cu^{1-s} \notag \\
&+\frac{2^{r}q_\cu^{1/2}}{\zeta_{\KK}^{\star}(2s)}\sum_{\substack{l\in \fb^{\star}\\ l \neq 0}}\tau_{1-2s}(l)K_{s-1/2}(\yy^{\star},l)e^{2\pi i Tr(l\xx^{\star})} \notag 
\end{align} where
\begin{align}
\tau_{s}(l)&=\NN(\fb\mathfrak{d}l)^{-s/2}\sum_{\substack{\mathfrak{q}  \text{ integral} \\
\mathfrak{q}: \fb\mathfrak{d}l }} \NN(\mathfrak{q})^{s} \notag
\end{align}
In addition, it satisfies the differential equation  $$\Delta \E(\zz,s)=(r_1+4r_2)
s(s-1)\E( \zz,s)$$ 
\end{coro}

\begin{pf}The Fourier expansion of $\E(\zz,s)$ at the cusp $\cu'$ follows by  dividing the Fourier expansion of $\E^{\star}(\zz,s)$ by $\zeta^{\star}_{\KK}(2s)$.
 Moreover, the differential equation $$\Delta(\NN(\yy)^s)=(r_1+4r_2)s(s-1)\NN(\yy)^s$$  and the  differential equation  of the Bessel function (together with its rapid decaying property) gives the differential equation for the analytic continuation of the Eisenstein series.
\end{pf}

\hfill

\begin{remark} 
The function 
\begin{equation}\label{funcionControladora}
\phi(s):=\frac{\zeta_{\KK}^{\star}(2s-1)}{\zeta_{\KK}^{\star}(2s)} \notag
\end{equation} which appears in the zero coefficient in different cusps, is of greatest importance as it governs many properties of the Eisenstein series.
\end{remark}

\section{Distribution of cusp cross sections}

\noindent
In this section we prove our statements. We first show the Rankin-Selberg unfolding trick for the Hilbert modular orbifold and the Maass-Selberg relation for the Eisenstein series $\E(\zz,s)$.
Distribution of the long closed horocycles and horospheres has been the subject of many works. For a dynamical point of view see \cite{V} and \cite{C}. 

\hfill

\noindent
Let $a\geq 0$ be an integer or infinity. We denote by $C_c^{a}(\MM)$  the set of complex valued functions defined on $\MM$ of class $C^a$ with compact support. Also $C(\MM)$ denotes the set of continuous complex valued functions on $\MM$. 
For the rest of the section we fix a  set of non equivalent  cusps $\{\lambda_i\}$ in order to simplify the arguments below. 
\hfill

\begin{definition}  Let $\X, \Y ,q$ be the local coordinates of $\zz$ at the cusp $\lambda_i$. 
For $f \in C(\MM)$ let $\widetilde{f}$ be a function in the variables $\X,\Y,q$  such that $\tilde{f}(\X,\Y,q)=f(\zz)$. The $\mathsf{measure}$ $\mathsf{uniformly}$ $\mathsf{supported}$ $\mathsf{at}$ $\mathsf{the}$ $\mathsf{cusp}$ $\mathsf{cross}$ $\mathsf{section}$ $\mathsf{at}$ $\mathsf{height}$ $q$ is defined by 
\begin{equation}\label{medida_cuspidal}
m_{i}(f,q)= \int_{-\frac{1}{2}}^{\frac{1}{2}}\cdots\int_{-\frac{1}{2}}^{\frac{1}{2}} \widetilde{f}(\X,\Y,q)d\X d\Y, \notag
\end{equation} 
where $d\X d\Y$ is the Lebesgue measure on $(\RR^{r_1} \times \CC^{r_2}) \times \RR^{r-1}$. 
\end{definition}
\hfill

\noi
For $f\in C_{c}^{a}(\MM)$ and each $i \in \{ 1,...,h\}$  consider the \emph{Mellin transform} of $m_i(f,q)q^{-1}$:
\begin{equation}\label{trasformadadeMellinmodular}
\mathcal{M}_{i}(f,s):=\int_{0}^{\infty}m_i(f,q)q^{s-1}\frac{dq}{q}. \notag
\end{equation}

\begin{prop} \label{convergenciadeM}
 For each $f\in C_c^a(\MM)$ and $i \in \{ 1,...,h\}$, the integral defining $\mathcal{M}_{i}(f,s)$ converges absolutely in the half-plane $\Re(s)>1$ and uniformly in  strips of the form $ \sigma_1 >\Re(s)> \sigma_{0}>1$. Therefore, it defines a  holomorphic function in the half plane $\Re(s)>1$.
\end{prop}
\begin{pf} 
  From the topology of the Hilbert modular orbifold it follows that
 a continuous function $f$  belongs to $C_c^{a}(\MM)$  if and only if  the corresponding  $\Gamma$-invariant function $\widetilde{f}$ in $\XX$ is of  class $C^a$ and there exists  a constant $\T>0$ such that if $\zz \in \XX$ satisfies $\mu(\zz,\lambda_i)>\T$, then $\widetilde{f}(\zz)=0$, for all $i=\{1,...,h\}$. 
Let $\norm{f}_{\infty}=\sup_{\substack{\zz\in \MM}}\{\abs{f(\zz)}\}$. Then, since $m_{i}(f,q) \leq \norm{f}_{\infty}$, for $\Re(s)>1$ we have:
\begin{equation}
\abs{\mathcal{M}_i(f,s)}\leq \norm{f}_{\infty}\left(\frac{\T^{\sigma-1}}{\sigma-1}\right), \qquad\text{where } \sigma=\Re(s) \notag
\end{equation} Therefore, we have absolute convergence in $\Re(s)>1$ and  uniform convergence in strips of the form $\sigma_1 >\Re(s) > \sigma_0>1$. 
\end{pf}

\hfill

\noi
The following theorem  is  the \emph{ Rankin-Selberg method} for the Hilbert modular orbifold. 
\hfill

\begin{prop}\label{MRS}  Let $f \in C_{c}^{0}(\MM)$ be a continuous function of compact support.
Then, if $\Re(s)>1$, we have
\begin{equation}\label{zetadeRS} 
\omega^{-1}\,2^{r_1-r_2}\,\R \, \sqrt{\D}\,\mathcal{M}_i(f,s)=\int_{\substack{\MM}}\E_{\cu_i}(\zz,s)\,f(\zz)\,d\upsilon(\zz) \notag
\end{equation}
\end{prop}
\begin{pf} Let $s$ be a point in the  half plane $\Re(s)>1$ and let $f$ be  of compact support. Then
\begin{align}
\int_{\substack{\MM}}\E_{\lambda_i}(\zz,s)\,f(\zz)\,d\upsilon(\zz)&=\sum_{\gamma\in\Gamma_{\lambda_i}\backslash \Gamma}
\int_{\MM}\mu(\gamma \zz,\lambda_i)^{s}\,f(\zz)\,d\upsilon(\zz) \notag
\end{align} 
Since the measure  $d\upsilon$ is invariant, changing $\MM$ for a  fundamental domain $F$,  we can see that the last series is equal to
$$
\sum_{\gamma\in\Gamma_{\lambda_i}\backslash \Gamma}\int_{\gamma(F)}\mu(\zz,\lambda_i)^{s}\,f(\zz)\,d\upsilon(\zz)
=\int_{\Gamma_{\lambda_i} \backslash \XX} \mu(\zz,\lambda_i)^s\,f(\zz)\,d\upsilon(\zz)$$
Then, if we use the local coordinates of $\zz$ at the cusp $\lambda_i$ to evaluate the integral, we have
\begin{align}
\int_{\substack{\Gamma_{\lambda_i}\backslash\XX}} \mu(z,\lambda_i)^s\,f(\zz)\,d\upsilon(\zz)&= \frac{2^{r_1}\,\R \,\sqrt{\D}}{2^{r_2}\,\omega}\int_{0}^{\infty}\int_{-\frac{1}{2}}^{\frac{1}{2}} \cdots\int_{-\frac{1}{2}}^{\frac{1}{2}}q^{s-2}\,\widetilde{f}(\X,\Y,q)\,d\X d\Y dq \notag \\
&= \omega^{-1}\,2^{r_1-r_2}\,\R \, \sqrt{\D}\,
\mathcal{M}_i(f,s) \notag
\end{align}
\end{pf}

\hfill

\noi
We shall need to estimate the Eisenstein series.
Let  $l_1(\KK)$ be the number given by  proposition \ref{independencia_cuspidal}, let  $\T>l_1(K)$ be a fixed large constant.  Define
$$
\E^{\T}(\zz,s)=\begin{cases} \E(\zz,s) & \text {if } q=\mu(\lambda,\zz) \leq \T \text{ for all cusps } \lambda \\

\E(\zz,s)-q^{s}-\phi(s)q^{1-s} & \text{if }  q=\mu(\lambda,\zz)>\T \text{ for some cusp } \lambda                  
                \end{cases}
$$

\hfill

\noi
The following theorem is the Maass-Selberg relation for the Eisenstein series of the Hilbert modular orbifold. We give a proof based on the unfolding trick  (cf. \cite[pp. 672]{Sel}).
\begin{prop}\label{Mass-Selberg} Let $s,s'$ be two complex numbers such that $s \neq s'$  $s+s'\neq 1$. Then
\begin{align}
\int_{\MM} \E^{\T}(\zz,s) \,\E^{\T}(\zz,s')\,d\upsilon(\zz)=&C\left(\frac{\T^{s+s'-1}-\phi(s)\phi(s')\T^{1-s-s'} }{s+s'-1}\right) \notag \\ 
&+C \left(\frac{\T^{s-s'}\phi(s')-\T^{s'-s}\phi(s)}{s-s'} \right) \notag
\end{align}
where $C=2^{r_1-r_2}\sqrt{\D} \, \R\, h \,\omega^{-1}$
\end{prop}
\begin{pf} 
First, from the Fourier expansion of  $\E^{\T}(\zz,s)$ at the cusps we see that for each cusp $\cu$ it decays quickly when $q_\cu \rightarrow \infty$ and uniformly in $s$ in vertical bands of finite width. Therefore   the left hand side of the above equation   is holomorphic in $s$  and $s'$ in the regular region of Eisenstein series.
Now, for a  real number $\T$ define 
 $$\delta_{\T}(q)=\begin{cases}
                   1 \text{ if } q > \T \\ 
                   0 \text{ if } q \leq \T
                  \end{cases}
$$ 
Let $F$ be a fundamental domain for the action of the modular group in $\XX$. 
For each $i=1,...,h$ and $\T$ large, we have an ``end'' of $F$:
 $$F_{\T}^i=\left\{ \zz \in F \ : \  q_i \geq \T \right\},$$
where $q_i=\mu(\zz,\lambda_i)$.
With this notation, if $\Re(s')>\Re(s)+1>2$, we have
\begin{align}
\int_{\MM} \E^{\T}(\zz,s)\, \E^{\T}(\zz,s')\,d\upsilon(\zz)& =\int_{F} \E(\zz,s) \,\E^{\T}(\zz,s') \,d\upsilon(\zz) \notag\\
&=
\int_{F} \E(\zz,s)\,\big( \E(\zz,s')-\sum_{i=1}^{h}\delta_{\T}(\q_i) \q_i^{s'}\big)\,d\upsilon(\zz)\notag\\ 
&-\phi(s')\,\sum_{i=1}^{h}\int_{F_{\T}^i} \E(\zz,s) \,\q_i^{1-s'}\,d\upsilon(\zz), \notag
\end{align} 
These last terms can be calculated separately as follows. First, recall that for $\Re(s)>1$, we have the bound $\abs{\E(\zz,s)}\leq \E(\zz,\sigma)$. Then 
\begin{align}
\int_{F} \E(\zz,s)\,\big( \E(\zz,s')-\sum_{\substack{i=1}}^{h}\delta_{\T}(\q_i)\,\q_i^{s'}\big)\, d\upsilon(\zz) =\notag\\ 
\int_{F} \E(\zz,s)\, \big(\E(\zz,s')-\sum_{\substack{i=1}}^{h}\,\q_i^{s'}\big)\,d\upsilon(\zz)+\sum_{\substack{i=1}}^{h}\int_{F - F_{\T}^i} \E(\zz,s) \q_i^{s'}\,d\upsilon(\zz) =\notag \\
\sum_{\substack{i=1}}^{h}\int_{F} \E(\zz,s) \sum_{\substack{\gamma \in \Gamma_{\lambda_i} \backslash (\Gamma-\Gamma_{\lambda_i})} }\mu(\lambda_i,\gamma\zz)^{s'}\,d\upsilon(\zz)+\sum_{\substack{i=1}}^{h}\int_{F-F_{\T}^{i}} \E(\zz,s) \q_i^{s'}\,d\upsilon(\zz)= \notag\\
\sum_{\substack{i=1}}^{h}\int_{F_{\lambda_i}-F} \E(\zz,s)\,\q_i^{s'}\,d\upsilon(\zz)+\sum_{\substack{i=1}}^{h}\int_{F-F_{\T}^i} \E(\zz,s) \q_i^{s'}\,d\upsilon(\zz) \notag,
\end{align}
where $F_{\lambda_i}$ is a fundamental domain for the action of $\Gamma_{\lambda_i}$ in $\XX$. 
Since the disjoint union of  $F_{\lambda_i}-F$ and $F-F_{\T}^i$  is equal to  
$$\left\{\zz:   \q_i  < \T \right\} \cap F_{\lambda_i},$$
 the last equation is equal to $2^{r_1-r_2}\sqrt{\D}\R  \omega^{-1}$ times
\begin{align}
\sum_{i=1}^{h} \int_{0}^{\T}\int_{-\frac{1}{2}}^{\frac{1}{2}} \cdots\int_{-\frac{1}{2}}^{\frac{1}{2}}\widetilde{\E}(\X_i ,\Y_i, \q_i)\q_i^{s'}\,d\X_i d\Y_i dq_i/q_i^{2}
&=\sum_{i=1}^{h}  \int_{0}^{\T}( q_i^{s}+ \phi(s)q_i^{1-s})\,q_i^{s'}\frac{d\q_i}{\q_i^{2}} \notag
\\
&=h \left(\frac{\T^{s+s'-1}}{s+s'-1}+\frac{\phi(s)\T^{s'-s}}{s'-s}\right) \notag
\end{align} 
On the other hand, the remaining term is equal to
\begin{align}
-\sum_{i=1}^{h}\phi(s')\int_{F_{T}^i} \E(\zz,s)\q_i^{1-s'} d\upsilon(\zz) & 
=-\frac{C}{h}\sum_{i=1}^{h}\phi(s')\int_{\T}^{\infty} (q_{i}^{s}+ \phi(s)q_{i}^{1-s})q_i^{1-s'}\frac{dq_i}{q_i^{2}} \notag \\
&=C\left(\frac{\phi(s')\T^{s-s'}}{s-s'}+\phi(s)\phi(s')\frac{\T^{1-s-s'}}{1-s-s'}\right) \notag
\end{align}  where $C=2^{r_1-r_2}\sqrt{\D} \, \R\, h \,\omega^{-1}$. This is valid for $\Re(s')>\Re(s)+1>2$.
Since both sides of the equation are meromorphic,  the identity follows for
 $s$ and $s'$ in the regular region of both sides.
\end{pf}

\hfill

\begin{remark}\label{cuota_alternativa}
 If we replace $\E(\zz,s)$ with the constant function $f=1$ in the procedure above,
we can evaluate (see \cite{Sel} and also \cite{DZ2} \cite{LW}) 
$$\int_{\MM_\T}\E(\zz,s')d\upsilon(\zz)=C\left(\frac{\T^{s'-1}}{s'-1}-\frac{\phi(s')\T^{-s'}}{s'}\right)
$$ for $\Re(s')>1$, where $C=2^{r_1-r_2}\sqrt{\D} \, \R\, h \,\omega^{-1}$. Since both sides are meromorphic functions, by analytic continuation,
this identity holds for any $s'$ in the regular region of both sides. 
Then, if we take the residue at $s'=1$ and let $\T$ tend to $\infty$, the volume of the Hilbert modular orbifold can be computed to be
$$vol(\MM)= 2^{-3r_2+1}\pi^{-n}\D^{3/2}\zeta_{\KK}(2).
$$
\end{remark}

\hfill

\noi
For each $q>0$  we 
have the measure $$m_q(f)=\frac{1}{h}\sum_{i=1}^{h}m_i(f,q)\qquad (f \in C_c^{0}(\MM))$$ and  the Mellin transform of $m_q(f)q^{-1}$:
\begin{equation}\label{trasformadadeMellinmodular2}
\mathcal{M}(f,s):=\int_{0}^{\infty}m_q(f)\,q^{s-1}\,\frac{dq}{q}=\frac{1}{h}\sum_{i}^{h}\mathcal{M}_i(f,s). \notag
\end{equation} 
Then, from  the  Rankin-Selberg theorem, for $f \in C_{c}^{0}(\MM)$ and $\Re(s)>1$, we have
$$
\omega^{-1}\,h\,2^{r_1-r_2}\,\R\,\sqrt{\D}\,\mathcal{M}(f,s)=\int_{\substack{\MM}}\E(\zz,s)\,f(\zz)\,d\upsilon(\zz) 
$$ Therefore, for a continuous function $f$ with compact support, the Mellin transform $\mathcal{M}(f,s)$  has the same properties of  $\E(\zz,s)$. That is to say, $\mathcal{M}_{f}(s)$ has a meromorphic continuation to the whole complex plane  that is regular for  $\Re(s)\geq \frac{1}{2}$ except, possibly, for a simple pole at $s=1$ with residue
\begin{equation}\label{resizetaRS}
\mathcal{R}es_{s=1} (\mathcal{M}_f(s))=\frac{1}{vol(\MM)} \int_{\MM}f(\zz)\,d\upsilon(\zz). \notag
\end{equation} 
From Lemma \ref{sii}, it follows that
the Riemann hypothesis for the field $\KK$ holds if and only if for all $f\in C_{c}^{0}(\MM)$ the function $\mathcal{M}_f(s)$ is regular for $\Re(s)>1/4$ except, possibly,  for a simple pole at $s=1$ with residue given as above.

\begin{remark}The modified function:
\begin{equation}
\mathcal{M}^{\star}_f(s)=\zeta^{\star}_{\KK}(2s)\,\mathcal{M}_f(s) \notag
\end{equation}   has a holomorphic continuation to the whole complex plane except, possibly, for simple poles at $s=0,1$ and satisfies the functional equation $$\mathcal{M}^{\star}_f(s)=\mathcal{M}^{\star}_f(1-s).$$
\end{remark}

\noi
Now, let us recall the following facts about the order of growth of $\zeta_{\KK}(s)$ along vertical lines. For each $\sigma$ we define a number $\mu(\sigma)$ as the lower bound of the numbers  $l\geq 0$ such that
\begin{equation}
 \zeta_{\KK}(\sigma+it)=\mathcal{O}(\abs{t}^{l}) \mbox{ as } \abs{t}\rightarrow \infty. \notag
\end{equation}
Then $\mu$ has the following properties \cite[p. 266]{SL}:
\begin{enumerate}[i)]
 \item    $\mu$ is continuous, non-increasing and never negative.
\item $\mu$  is convex downwards in the sense that the curve $y=\mu(\sigma)$ 
 has no points above the chord joining any two of its points.
\item  $\mu(\sigma)=0$   if  $\sigma \geq 1$  and  $\mu(\sigma)=n(\frac{1}{2}-\sigma)$ if  $\sigma  \leq 0$.
\end{enumerate}

\hfill

\begin{lemma}\label{cotas_principal}
Let $f$ be a  differentiable function of compact support in  $\MM$. Then,  for each  $1/2<\sigma_0<1$, there exists $t_0$ such that
$$
(r_1+4r_2)\abs{\M_f(s)}\leq \frac{\beta_f(\sigma_0)\,t^{\epsilon} }{\abs{s(s-1)}} \quad \text{ for } t>t_0, \ \epsilon>0
$$
and $\sigma_0\leq\Re(s)\leq 2$,
where  $\beta_f(\sigma_0)$ is a constant depending on $\sigma_0$ and the second  derivatives of $f$.
 \end{lemma}

\begin{pf} From proposition (\ref{MRS}) we have
$$\mathcal{M}(f,s)=\int_{\substack{\MM}}\E(\zz,s)\,f(\zz)\,d\upsilon(\zz).$$
Since $ \Delta \E(\zz,s) =(r_1+4r_2)s(s-1)\E(\zz,s)$ and $\Delta$ is a  symmetric operator, then
\begin{align}
(r_1+4r_2)s(s-1)\mathcal{M}(f,s)&=\int_{\substack{\MM}}\Delta\E(\zz,s)\,f(\zz)\,d\upsilon(\zz) \notag\\
 &=\int_{\substack{\MM}}\E(\zz,s)\Delta f(\zz)d\omega(\zz) \notag
\end{align}
Now, let $\T>0$ be such that $ supp f \subset \{\, q_j \leq \T \,\}$. Thus,
\begin{align} 
(r_1+4r_2)s(s-1)\mathcal{M}(f,s) &=\int_{\substack{\MM}}\E^{\T}(\zz,s)\,\Delta f(\zz)\,d\upsilon(\zz) \notag
\end{align}
Therefore, the Cauchy-Schwarz inequality implies
$$(r_1+4r_2)\abs{s(s-1)} \abs{\mathcal{M}(f,s)} \leq \norm{\Delta f}_2 \, \norm{\E^{\T}(\cdot,s)}_2.
$$
Now, from the Mass-Selberg relations, the $L^2$ norm of the Eisenstein series  is given by 
\begin{align}
\int_{\MM} \abs{\E^{\T}(\zz,s)}^{2}d\upsilon(\zz)&=C\left(\frac{\T^{2\sigma-1}-\abs{\phi(s)}^2\T^{1-2\sigma} }{2\sigma-1} \right) \notag \\
&+ C\left(\frac{\T^{2it}\phi(\overline{s})-\T^{-2it}\phi(s)}{2it} \right) \notag 
\end{align} where $s=\sigma+it, t \neq 0 , \sigma \neq \frac{1}{2}$  and  $C=2^{r_1-r_2}\, \R\,\sqrt{\D}\,h\,\omega^{-1}$.
Therefore, we can bound the $L^2$ norm of the truncated Eisenstein series in vertical bands of finite width essentially by   $\phi(s)$, except for $\sigma=1/2$ and $t=0$.
 From Stirling's formula (cf. \cite[Ch. 13]{WW}), it follows that 
$$
\Gamma(\sigma+it)\sim\sqrt{2\pi} \abs{t}^{\sigma-1/2}e^{-\pi \abs{t}/2}
$$ as $t \rightarrow \pm \infty$ uniformly in vertical bands of finite width. Thus,
$$\frac{2^{r_2}\pi^{\frac{n}{2}}\Gamma(s-\frac{1}{2})^{r_1}\Gamma(2s-1)^{r_2}}{\sqrt{\mathrm{D}}\Gamma(s)^{r_1}\Gamma(2s)^{r_2}}
 \sim \pi^{\frac{n}{2}} \D^{-\frac{1}{2} }\abs{t}^{-n/2}
$$ as $t \rightarrow \pm \infty$ uniformly in $1/2 \leq \Re(s) \leq 2$. On the other hand, if $l>n/2=\mu(0)$, then
$$
\zeta_{\KK}(2s-1)=\mathcal{O}(t^{l}),
$$
uniformly in $1/2 \leq \Re(s) \leq 2$. Furthermore, if $ \sigma_0<\Re(s)$, thus
$$
\frac{1}{\zeta_{\KK}(2s)}=\sum_{\substack{ \fa \subset \OO}}\frac{\mu(\fa)}{\NN(\fa) ^{2s}}
$$ where $\mu(\fa)$ is the M\"obius function of the field $\KK$.
Therefore, for all $\epsilon >0$,  we have
\begin{equation}
 \abs{\phi(s)}=\mathcal{O}(\abs{t}^{\epsilon}) \notag
\end{equation} as $\abs{t} \rightarrow \infty$ uniformly in $\sigma_0\leq \Re(s)\leq 2$, which proves the claim.
\end{pf}

\begin{p1}
Let $f \in C_{c}^{\infty}(\MM)$. Then, from the Ranking-Selberg representation, the Mellin transform $\mathcal{M}(f, s)$ of $m_q(f)q^{-1}$  is holomorphic for $\Re(s)>\frac{1}{2} $ 
except, possibly, for a simple pole at $s=1$ with residue
$$
m(f)=\frac{1}{vol(\MM)} \int_{\MM} f(\zz)d\upsilon(\zz).
$$
From proposition \ref{convergenciadeM} the Mellin inversion formula applies, and we have
\begin{equation}\label{MI}
m_q(f)=\frac{1}{2\pi i}\int_{b-i\infty}^{b+i\infty}\mathcal{M}_f(s)q^{1-s}ds. 
\end{equation} for any real number $b>1$.
Now, let $0< \epsilon <\frac{1}{2}$.  Then, by the estimates of Lemma \ref{cotas_principal}, we  can shift the path of integration in equation (\ref{MI}) to the line
 $\sigma=\frac{1}{2}+\epsilon$ to get
$$
m_q(f)= m(f) +\frac{1}{2\pi}\int_{-\infty}^{\infty}\mathcal{M}_f(\frac{1}{2}+\epsilon+it)q^{\frac{1}{2} +\epsilon} q^{-it}dt.
$$
Now, since $\mathcal{M}_f(\frac{1}{2}+\epsilon+it)$ is integrable (w.r.t. $dt$), the Riemann-Lebesgue theorem implies
$$
\lim_{q\rightarrow 0}\abs{\int_{-\infty}^{\infty}\mathcal{M}_f(\frac{1}{2}+\epsilon+it)q^{it}dt}=0
$$ 
Therefore,
$$
m_q(f)=m(f)+o(q^{1/2-\epsilon}) \qquad (q \rightarrow 0).
$$ 
\end{p1}

\hfill

\noi
Now we show the relation with the Riemann hypothesis of the Dedekind zeta function.
First, we state a classical consequence of the Riemann hypothesis (cf. \cite[p. 267]{SL}). 
\begin{prop}\label{implicaciondeHR} 
If the Riemann hypothesis for the Dedekind zeta function of $\KK$ is true,
then we have the following estimates
\begin{gather}
\text{For } \epsilon>0 \text{ and }\sigma > \frac{1}{2}:  \notag\\
- \epsilon \log t < \log \abs{\zeta_{K}(s)} < \epsilon \log t; \ s=\sigma+it, \ t \geq t_0(\epsilon), \notag
\end{gather}
 that is to say,1
\begin{equation}\label{Littlewood}
\begin{cases}
\zeta_{\KK}(s)=\mathcal{O}(t^{\epsilon}) & \\
&  \text{for every } \epsilon>0, \ s=\sigma+it, \sigma>\frac{1}{2} \text{ as }\abs{t} \rightarrow \infty.\\
 \frac{1}{\zeta_{\KK}(s)}=\mathcal{O}(t^{\epsilon})&
\end{cases} \notag
\end{equation}
\end{prop}
\hfill

\begin{lemma}\label{cotas_principal_HR}  If the Riemann hypothesis for the Dedekind zeta function holds, then for each $1/4<\sigma_0<1/2$, there exists $t_0>0$ such that
$$
\abs{\M_f(s)}\leq \frac{\beta_f(\sigma_0)}{\abs{s(s-1)}}, \qquad \text{ for } t>t_0,
$$
for all $s$ with $\sigma_0\leq\Re(s)\leq 2$ and $\Re(s)\neq 1/2$. Here $\beta_f(\sigma_0)$ is a constant depending essentially on $\sigma_0$ and a finite number of derivatives of $f$. 
\end{lemma}
\begin{pf} First, we estimate $\phi(s)=\zeta^{\star}_{\KK}(2s-1)/\zeta_{\KK}^{\star}(2s)$ in the region $\sigma_0\leq\Re(s) \leq 2$, under the assumption of the Riemann hypothesis for the Dedekind zeta function $\zeta_{\KK}(s)$. From proposition \ref{implicaciondeHR}, for every $\epsilon>0$,  $\zeta_{\KK}(2s)^{-1}=\mathcal{O}(t^{\epsilon})$, uniformly in $\sigma_0\leq\Re(s) \leq 2$. Alike, $\zeta(2s-1)=\mathcal{O}(t^l)$ for $l\geq n=\mu(-1/2)$, uniformly in $\sigma_0\leq\Re(s) \leq 2$.
Therefore, Stirling formula implies that, for every $\epsilon>0$, $\phi(s)=\mathcal{O}(t^{n/2+\epsilon})$, uniformly in $\sigma_0\leq\Re(s) \leq 2$. Now  we can use the process of Lemma \ref{cotas_principal}  to see that a sufficient degree of
derivatives of $f$ ensures that, 
$$\abs{\mathcal{M}_f (s)} =\mathcal{O}\big(\abs{s(s-1)}^{-1} \big) \qquad(\abs{
s}\rightarrow \infty)$$ 
for $\Re(s)\neq 1/2$ and uniformly in $\sigma_0\leq\Re(s) \leq 2$. This proves the assertion.
\end{pf}

\hfill

\begin{p2} 
Suppose that for all  $f\in C_c^{\infty}(\MM)$, we have the following bound:
$$m_{q}(f)=m(f)+\mathcal{O}(q^{3/4-\epsilon}) \qquad(q\rightarrow 0)
$$ for all $0<\epsilon<3/4$
and write $m_q(f)=m(f)+k(q)$. Let   $\T$ be sufficiently large and such that $m_q(f)=0$ for $q>\T$.  Then,
\begin{align}
\mathcal{M}_{f}(s)&=\int_{0}^{\T}m_{q}(f)q^{s-2}dq \notag\\
&=\int_{0}^{\T}(m(f)+k(q))q^{s-2}dq \notag \\
&=\frac{m(f)\T^{s-1}}{s-1}+~\int_{0}^{\T}k(q)q^{s-2}dq \notag
\end{align} 
Since $k(q)=\mathcal{O}(q^{\frac{3}{4}-\epsilon})$, the last integral 
 converges absolutely and uniformly in the half-plane $\Re(s)>\frac{1}{4}+\epsilon$, so it defines a holomorphic function in that half-plane. Therefore, 
$\mathcal{M}_f(s)$ is a holomorphic function in the region $\Re(s)>\frac{1}{4}+\epsilon$ except, possibly, for a pole  at $s=1$ with residue $m(f)$. Thus, the Riemann hypothesis for  the Dedekind zeta function is true. On the other hand,  suppose the Riemann hypothesis for the Dedekind zeta function of  $\KK$ holds. Then,  $\zeta_{\KK}(2s)$ does dot vanish  for $\Re(s)>1/4$ and
$\mathcal{M}_f(s)$ is holomorphic for $\Re(s)>1/4$ except, possibly, for a simple pole at $s=1$ with residue $m(f)$. From Lemma \ref{cotas_principal},  the  integral of $\mathcal{M}_f(s)q^{1-s}$ exists over the boundary of the band $\frac{1}{4}+\epsilon \leq \sigma \leq 2$, for all $0<\epsilon< 1/4$. Therefore, the Mellin inversion formula  implies
$$
m_q(f)=\mathcal{R}es_{s=1}(\mathcal{M}_f(s))+\frac{1}{2\pi}\int_{-\infty}^{\infty}\mathcal{M}_f(\frac{1}{4}+\epsilon +it)q^{-it}q^{\frac{3}{4}-\epsilon}dt.
$$ Again, by the Riemann-Lebesgue Theorem:
$$
m_q(f)=m(f)+o(q^{\frac{3}{4}-\epsilon})  \qquad (q\rightarrow 0).
$$
\end{p2}

\section*{Acknowledgements} 
\noi
The author would like to thank  Dr. Alberto Verjovsky for suggesting the research of Zagier's distribution theorem on the Hilbert  modular orbifold
 and for all his great  explanations on this remarkable subject.

\bibliographystyle{model1b-num-names}

\bibliography{references}

\begin{thebibliography}{23}
\expandafter\ifx\csname natexlab\endcsname\relax\def\natexlab#1{#1}\fi
\providecommand{\bibinfo}[2]{#2}
\ifx\xfnm\relax \def\xfnm[#1]{\unskip,\space#1}\fi
\bibitem[{Bump(1998)}]{DB}
\bibinfo{author}{D.~Bump}, \bibinfo{title}{Automorphic Forms and
  Representations}, \bibinfo{publisher}{Cambridge University Press},
  \bibinfo{address}{Cambridge}, \bibinfo{year}{1998}.
\bibitem[{Cohen and Sarnak(1980)}]{CP}
\bibinfo{author}{P.~Cohen}, \bibinfo{author}{P.~Sarnak},
  \bibinfo{title}{Eisenstein series for hyperbolic manifolds},
  \bibinfo{year}{1980}. \bibinfo{note}{Notes with P. Cohen on the trace formula
  Chapter 6 and 7, Private notes, Princeton
  http://www.math.princeton.edu/sarnak/}.
\bibitem[{Cosentino(1999)}]{C}
\bibinfo{author}{S.~Cosentino}, \bibinfo{title}{Equidistribution of parabolic
  fixed points in the limit set of {K}leinian groups},
  \bibinfo{journal}{Ergodic Theory and Dynamical Systems}
  (\bibinfo{year}{1999}) \bibinfo{pages}{1437--1484}.
\bibitem[{Efrat(1987)}]{Ef}
\bibinfo{author}{I.~Efrat}, \bibinfo{title}{The Selberg trace formula for
  $PSL(2,\mathbb{\RR})^n$}, \bibinfo{publisher}{Memoirs of the American
  Mathematical Society}, \bibinfo{address}{Providence Rhode Island USA},
  \bibinfo{year}{1987}.
\bibitem[{Elstrodt et~al.(1997)Elstrodt, Grunewald and Mennicke}]{ELS}
\bibinfo{author}{J.~Elstrodt}, \bibinfo{author}{F.~Grunewald},
  \bibinfo{author}{J.~Mennicke}, \bibinfo{title}{Groups Acting on Hyperbolic
  Space}, \bibinfo{publisher}{Springer-Verlag}, \bibinfo{year}{1997}.
\bibitem[{Estala(2012)}]{Sam}
\bibinfo{author}{S.~Estala}, \bibinfo{title}{Distribuci\'on de secciones
  cuspidales en la variedad modular de Hilbert}, Ph.D. thesis, Universidad
  Nacional Aut\'onoma de M\'exico, \bibinfo{year}{2012}. \bibinfo{note}{To
  appear}.
\bibitem[{Freitag(1990)}]{F}
\bibinfo{author}{E.~Freitag}, \bibinfo{title}{Hilbert Modular forms},
  \bibinfo{publisher}{Springer-Verlag}, \bibinfo{address}{Berlin-New York},
  \bibinfo{year}{1990}.
\bibitem[{van~der Geer(1988)}]{GV}
\bibinfo{author}{G.~van~der Geer}, \bibinfo{title}{Hilbert Modular Surfaces},
  \bibinfo{publisher}{Springer-Verlag}, \bibinfo{address}{Berlin-New York},
  \bibinfo{year}{1988}.
\bibitem[{Hejhal(1983)}]{He}
\bibinfo{author}{D.A. Hejhal}, \bibinfo{title}{The Selberg Trace Formula for
  $PSL(2,\mathbb{R})$}, volume~\bibinfo{volume}{II},
  \bibinfo{publisher}{Springer- Verlag}, \bibinfo{address}{Berlin-New York},
  \bibinfo{year}{1983}.
\bibitem[{Jorgenson and Lang(1999)}]{L2}
\bibinfo{author}{J.~Jorgenson}, \bibinfo{author}{S.~Lang},
  \bibinfo{title}{Hilbert-{A}sai {E}isenstein series, regularized products, and
  heat kernels}, \bibinfo{journal}{Nagoya Math. J.}  (\bibinfo{year}{1999})
  \bibinfo{pages}{155--188}.
\bibitem[{Kubota(1973)}]{Ku}
\bibinfo{author}{T.~Kubota}, \bibinfo{title}{Elementary Theory of Eisenstein
  Series}, \bibinfo{publisher}{Kodansha LTD Japan Tokyo, Halsted press USA},
  \bibinfo{year}{1973}.
\bibitem[{Lang(1994)}]{SL}
\bibinfo{author}{S.~Lang}, \bibinfo{title}{Algebraic Number Theory},
  \bibinfo{publisher}{Springer-Verlag}, \bibinfo{address}{Berlin-New York},
  \bibinfo{year}{1994}.
\bibitem[{Neukirch(1999)}]{JN}
\bibinfo{author}{J.~Neukirch}, \bibinfo{title}{Algebraic Number Theory},
  \bibinfo{publisher}{Springer-Verlag}, \bibinfo{address}{Berlin-New York},
  \bibinfo{year}{1999}.
\bibitem[{Sarnak(1981)}]{Sa}
\bibinfo{author}{P.~Sarnak}, \bibinfo{title}{Asymptotic behavior of periodic
  orbits of the horocycle flow and {E}isenstein series},
  \bibinfo{journal}{Comm. on Pure and Applied Math.}  (\bibinfo{year}{1981})
  \bibinfo{pages}{719--739}.
\bibitem[{Selberg(1991)}]{Sel}
\bibinfo{author}{A.~Selberg}, \bibinfo{title}{Collected papers},
  volume~\bibinfo{volume}{I}, \bibinfo{publisher}{Springer-Verlag},
  \bibinfo{year}{1991}.
\bibitem[{Siegel(1961)}]{S}
\bibinfo{author}{C.~Siegel}, \bibinfo{title}{Lectures notes on Advanced
  Analytic Number Theory}, \bibinfo{publisher}{Tata Institute},
  \bibinfo{address}{Bombay}, \bibinfo{year}{1961}.
\bibitem[{Sorensen(2002)}]{CM}
\bibinfo{author}{C.M. Sorensen}, \bibinfo{title}{Fourier expansion of
  {E}isenstein series on the {H}ilbert modular group and {H}ilbert class
  fields}, \bibinfo{journal}{Transactions of the American Mathematical Society}
   (\bibinfo{year}{2002}) \bibinfo{pages}{4847--4869}.
\bibitem[{Thurston(2002)}]{Th}
\bibinfo{author}{W.P. Thurston}, \bibinfo{title}{The geometry and topology of
  three-manifolds}, \bibinfo{year}{2002}. \bibinfo{note}{Lecture Notes,
  Princeton University, http://library.msri.org/books/gt3m/}.
\bibitem[{Verjovsky(1992)}]{V}
\bibinfo{author}{A.~Verjovsky}, \bibinfo{title}{Arithmetic, geometry and
  dynamics in the modular orbifold}, \bibinfo{journal}{Dynamical Systems,
  (Santiago de Chile 1990)(Pitman Series 285), R. Bamon, R. Labarca. J.
  Lewowicz, J. Palis, Longman, Essex, UK}  (\bibinfo{year}{1992})
  \bibinfo{pages}{263--298}.
\bibitem[{Weng(2006)}]{LW}
\bibinfo{author}{L.~Weng}, \bibinfo{title}{A rank two zeta and its zeros},
  \bibinfo{journal}{J. Ramanujan Math. Soc. 21, No.3}  (\bibinfo{year}{2006})
  \bibinfo{pages}{1--62}.
\bibitem[{Whittaker and Watson(1927)}]{WW}
\bibinfo{author}{E.~Whittaker}, \bibinfo{author}{G.~Watson}, \bibinfo{title}{A
  Course of Modern Analysis}, \bibinfo{publisher}{Cambridge University Press},
  \bibinfo{address}{Cambridge}, \bibinfo{edition}{fourth} edition,
  \bibinfo{year}{1927}.
\bibitem[{Zagier(1979)}]{DZ}
\bibinfo{author}{D.~Zagier}, \bibinfo{title}{Eisenstein series and the
  {R}iemann zeta function}, \bibinfo{journal}{Automorphic forms, Representation
  theory and Arithmetic, Tata Institute of Fundamental Research, Bombay}
  (\bibinfo{year}{1979}).
\bibitem[{Zagier(1981)}]{DZ2}
\bibinfo{author}{D.~Zagier}, \bibinfo{title}{The {R}ankin-{S}elberg method for
  automorphic functions which are not of rapid decay}, \bibinfo{journal}{J.
  Fac. Sci. Univ. Tokyo}  (\bibinfo{year}{1981}).

\end{thebibliography}

\end{document}